\documentclass[12pt,aos]{imsart}

\usepackage{amsthm,amsmath}

\usepackage{amsfonts, amssymb}
\usepackage{mathrsfs}
\usepackage{epsfig}
\usepackage{graphicx}
\usepackage{amsmath}
\usepackage{color}
\newtheorem{prop}{Proposition}
\newtheorem{thm}{Theorem}
\newtheorem{cor}{Corollary}
\newtheorem{lem}{Lemma}
\newtheorem{ex}{Example}
\startlocaldefs

\def\baselinestretch{1.5}
\def\B{\mathbb B}
\def\M{\mathbb M}

\def\R{\mathbb R}
\def\W{\mathbb W}
\def\X{\mathbb X}
\def\Y{\mathbb Y}
\def\Z{\mathbb Z}
\def\qed{\quad $\Box$}
\def\ed{\end{document}}

\def\be{\begin{equation}}
\def\ee{\end{equation}}
\def\ed{\end{document}}

\def\<{\langle}
\def\>{\rangle}
\def\eqalign#1{\null\,\vcenter{\openup\jot\ialign
              {\strut\hfil$\displaystyle{##}$&$\displaystyle{{}##}$
               \hfil\crcr#1\crcr}}\,}

\endlocaldefs

\begin{document}

\begin{frontmatter}

\title{Estimating a Monotone Trend}
\runtitle{Estimating a Monotone Trend}


\author{\fnms{Ou} \snm{Zhao}\ead[label=e1]{ou.zhao@yale.edu}}
\address{Division of Biostatistics \\
Yale University \\
School of Medicine \\
New Haven, Connecticut \\
USA \\
\printead{e1}}
\and
\author{\fnms{Michael} \snm{Woodroofe}\ead[label=e2]{michaelw@umich.edu}}
\address{Department of Statistics \\
University of Michigan \\
275 West Hall \\
Ann Arbor, MI \\
USA \\
\printead{e2}}
\affiliation{Yale University and University of Michigan}

\runauthor{O. Zhao and M. Woodroofe}

\begin{abstract}
Motivated by global warming issues, we consider a time series that consists of a
nondecreasing trend observed
with stationary fluctuations, nonparametric estimation of the trend
under monotonicity assumption is considered. The rescaled  isotonic
estimators at an interior point are shown to converge to Chernoff's distribution
under minimal conditions on the stationary errors. Since the isotonic
estimators suffer from the spiking problem at the end point, two modifications
are proposed. The estimation errors for both estimators of the boundary point
are shown to have interesting limiting distributions. Approximation
accuracies are  assessed through simulations. One highlight of our treatment
is the proof of the weak convergence results which involve several recent
techniques developed in the study of conditional central limit questions. These weak
convergences can be shown to hold conditionally given the starting values.

\end{abstract}



\end{frontmatter}

\section{Introduction}\label{sect:intro}

	Consider a time series that consists of a nondecreasing trend observed with stationary fluctuations,  say
$$
	y_k = \mu_k + X_k,\quad k = 1,2,\ldots
$$
where $-\infty < \mu_1 \le \mu_2 \leq \cdots$ and $\ldots X_{-1},X_0,X_1,\ldots$ is a strictly stationary sequence with mean $0$ and finite variance.  The global temperature anomalies in Example 1 provide a particular example. If a segment of the series is observed, say $y_1,\ldots,y_n$, then isotonic methods suggest themselves for estimating the $\mu_k$ nonparametrically.  The isotonic estimators may be described as
\be\label{eq:minmax}
	\tilde{\mu}_k = \max_{ i\le k} \min_{k \le j \le n} {y_i+\cdots+y_j\over j-i+1}.
\ee
Alternatively, letting $\lfloor x\rfloor$ denote the greatest integer that is less than or equal to $x \in \R$, $Y_n$ the cumulative sum diagram, 
$$
	Y_n(t) = {y_1+\cdots+y_{\lfloor nt\rfloor} \over n},
$$
and $\tilde{Y}_n$ its greatest convex minorant, $\tilde\mu_k = \tilde{Y}_n^{'}(k/n)$, the left hand derivative of $\tilde{Y}_n$ evaluated at $t = k/n$.  See Chapter 1 of \cite{RWD88} for background on isotonic estimation.

\begin{ex}\label{ex:glbltmp} 
	{\rm Annual global temperature anomalies from 1850-2000 are shown in  Figure \ref{fig:glbltmp} with the isotonic estimator of trend superimposed as a step function.

\begin{figure}[htbp] 
   \centering 
   \includegraphics[width=3.7in,angle=270]{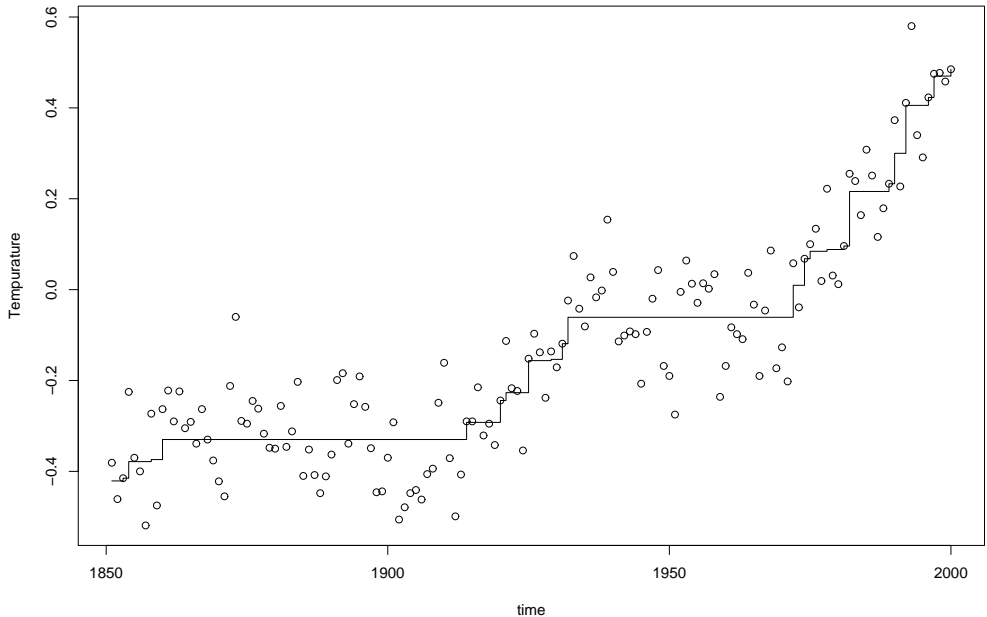} 
   \caption{Global Temperature Anomalies}
   \label{fig:glbltmp}
\end{figure}

}
\end{ex}

	With the global warming data, there is special interest in estimating $\mu_n$, the current temperature anomaly, and there isotonic methods encounter the {\it spiking problem}, described in Section 7.2 of \cite{RWD88} for the closely related problem of estimating a monotone density.  We consider two methods for correcting this problem, the penalized estimators of \cite{WS93} and the method of \cite{KL06}, both introduced for monotone densities.  The former estimates $\mu_n$ by
$$
	\hat\mu_{p,n} = \max_{i\le n} {y_i+\cdots+y_n\over n-i+1+\lambda_n},
$$
where $\lambda_n > 0$ is smoothing parameter, and the latter by $\hat\mu_{b,n} = \tilde\mu_{m_n}$,  where $m_n < n$ is another smoothing parameter.  

	The main results of this paper obtain the asymptotic distributions of estimation errors, properly normalized, for the estimators described above.  One of these results is well known for monotone regression with i.i.d. errors, and analogues of the others are known for monotone density estimation.   Interest here is in extending these results to allow for dependence.  Others have been interested in this  question recently---notably Anevski and H\" ossjer \cite{AH06}.  Our results go beyond theirs in several ways.  We consider the boundary case, estimating $\mu_n$; our results hold conditionally given the starting values; and our conditions are weaker.  Instead of the strong mixing condition, called (A9) in \cite{AH06}, we use the condition (\ref{eq:MW}) below, introduced in \cite{MW00} and further developed in \cite{PU05}.  One objective of this paper is to show by example how recent results on the central limit question for sums of stationary processes can be used to weaken mixing conditions in statistical applications.

	The main results are stated and proved in Section \ref{sect:ad} and then illustrated by simulations in Section {\ref{sect:sims}.  Section \ref{sect:prelim} contains some background material.

\section{Preliminaries}\label{sect:prelim} \hspace{1cm}

	{\it A maximal inequality and conditional convergence}.\quad  The main results of \cite{PU05} are an important technical tool.  To state them, let $\ldots X_{-1},X_0,X_1,\ldots$ be a strictly stationary sequence with mean $0$ and finite variance, as above; let $S_0=0$, $S_n = X_1+\cdots+X_n$, ${\cal F}_n = \sigma\{\ldots,X_{n-1},X_n\}$, and 
$$
	\B_n(t) = {1\over\sqrt{n}}S_{\lfloor nt\rfloor}
$$
for $0 \le t \le 1$; and let $\B$ denote a standard Brownian motion.  Both $\B$ and  $\B_n$ are regarded as random elements with values in $D[0,1]$, 
 endowed with the Skorohod topology, \cite{B68}, Chapter 3.  Let $\Vert\cdot\Vert$ denote the norm in $L^2(P)$, $\Vert Y\Vert = \sqrt{E(Y^2)}$.  It is 
 shown in \cite{PU05} that if 
\be\label{eq:MW}
	\sum_{n=1}^{\infty} n^{-{3\over 2}}\Vert E(S_n|{\cal F}_0)\Vert  < \infty,
\ee
then
$$
 \Gamma :=\sum_{k=0}^{\infty} 2^{-{1\over 2} k} \Vert E(S_{2^k}\vert \mathcal{F}_0) \Vert <\infty ,
$$
and
\be\label{eq:PU}
	E \left[\max_{k\le n} S_k^{2} \right] \le 6 \left[ E(X_1^{2})+\Gamma \right] n,
\ee
\be\label{eq:sigma}
	\sigma^2 = \lim_{n\to\infty} {1\over n}E(S_n^{2})
\ee
exists, and $\B_n$ converges in distribution to $\sigma \B$.  In fact, a stronger conclusion is possible.  It will be shown that the conditional distributions of $\B_n$ given $\mathcal{F}_0$ converge in probability to the distribution of $\sigma\B$.

	Properties of weak convergence---for example, the continuous mapping theorem and  Slutzky's theorem, extend easily to the convergence of conditional distributions.  We illustrate with Slutzky's theorem \cite{B68}. Let $({\cal X},d)$ denote a complete separable metric space, and let $\rho$ be a metric that metrizes weak convergence of probability distributions on the Borel sets of ${\cal X}$, for example the metric (\ref{eq:rho}) below.  Next, let $\X_n, \Y_n,
	 n=1,2,\ldots$, be random elements assuming values in ${\cal X}$; suppose that $\X_n$ and $\Y_n$ are defined on the same probability space $(\Omega_n,{\cal A}_n,P_n)$ say; let ${\cal A}_n^{o} \subseteq {\cal A}_n$ be sub sigma algebras; and let $\mu_n$ and $\nu_n$ be regular conditional distributions for $\X_n$ and $\Y_n$ given ${\cal A}_n^{o}$. {\it If $\rho(\mu,\mu_n) \to 0$ in probability and $d(\X_n,\Y_n) \to 0$ in probability, then $\rho(\mu,\nu_n) \to 0$ in probability.}  The assertion can be easily proved from the usual statement of Slutzky's theorem, for example, \cite{B68}, page 25, by considering subsequences which converge to $\infty$ so rapidly that $\rho(\mu,\mu_n) \to 0$ and $d(\X_n,\Y_n) \to 0\ w.p.1$ along the subsequence. 

	There is a convenient choice of $\rho$.  Write $\Vert g\Vert_{{\rm Lip}} = \sup_x |g(x)| + \sup_{x\ne y} |g(x)-g(y)|/d(x,y)$ for bounded Lipschitz continuous functions $g$ and let
\be\label{eq:rho}
	\rho(\mu,\nu) = \sup_{\Vert g\Vert_{\rm Lip} \le 1} \left| \int_{\cal X} gd\mu - \int_{\cal X} gd\nu\right|
\ee
for probability distributions $\mu$ and $\nu$ on the Borel sets of ${\cal X}$. Then $\rho$ metrizes convergence in distribution (\cite{D89}, Theorem 11.3.3).   Here is a useful feature of $\rho$.  Let ${\cal A}_{n,1} \subseteq {\cal A}_{n,2}$ be sub sigma algebras of ${\cal A}_n$ and let $\mu_{n,1}$ and $\mu_{n,2}$ be regular conditional distributions for $\X_n$ given ${\cal A}_{n,1}$ and ${\cal A}_{n,2}$.  Then $\rho(\mu,\mu_{n,1}) \le E[\rho(\mu,\mu_{n,2} ) \vert {\cal A}_{n,1}]$ and, therefore,
\be\label{eq:rhoprp}
		E[\rho(\mu,\mu_{n,1} ) ] \le E[\rho(\mu,\mu_{n,2} )].
\ee 

	One more bit of preparation: If $\ldots X_{-1},X_0,X_1,\ldots$ is any stationary sequence for which $E(X_n^{2}) < \infty$, then $X_n^{2}/n \to 0\ w.p.1$ by an easy application of the Borel-Cantelli lemmas, or simply the ergodic theorem \cite{P89}, page 30;  thus,  $X_n/\sqrt{n} \to 0 $ with probability one.  

	If $\gamma > 0$, and $m \ge 0$ is an integer, let
\be\label{eq:wab}
	\X_{m,\gamma} (t) = {S_{m+\lfloor\gamma t\rfloor}-S_{m}\over\sqrt{\gamma}}
\ee
for $ t\geq -m/ \gamma$; and let $\X_{m,\gamma}^{a,b}=\X_{m,\gamma} \vert [a,b]$  denote the restriction of $X_{m,\gamma}$ to an interval $[a,b]$. Thus $\B_n = \X_{0,n}^{0,1}$.  Let $\W$ denote a standard two-sided Brownian motion. Both $\X_{m,\gamma}^{a,b}$ and $\W^{a,b}=\W \vert [a,b]$ are regarded as 
elements of $D[a, b]$.

\begin{prop}\label{prop:cc}
	Suppose that (\ref{eq:MW}) holds; let $m_n \ge 0$ be integers; let $0 < \gamma_n \to \infty$; and let $-\infty < a < b < \infty$.  If either $a \ge 0$ or $\gamma_n/m_n \to 0$, then the conditional distribution of $\X_{m_n,\gamma_n}^{a,b}$ given ${\cal F}_0$ converges in probability to the distribution of 
	$\sigma \W^{a,b}$. 
\end{prop}
	
	{\sc Proof}.\quad For fixed  $a$ and $b$,  write $\X_n = \X_{m_n,\gamma_n}^{a,b}$; let $\Phi_n$ denote a regular conditional distribution for $\X_n$ given ${\cal F}_0$ and  $\Phi$ the distribution of $\sigma\W^{a,b}$.  Then it is necessary to show that  $\rho[\Phi,\Phi_n] \to  0$ in probability.
	
	If $a \ge 0$, then it suffices to consider the case $a = 0$, since then the convergence of $\X_{m_n,\gamma_n}^{0,b}$ implies that of $\X_{m_n,\gamma_n}^{a,b}$.  It also suffices to consider the case $m_n = 0$.  To see why, suppose that the result is known for $m_n = 0$ and let $\Phi_{n}^{o}$ be a regular conditional distribution for $\X_{0,\gamma_n}^{0,b}$ given ${\cal F}_0$, so that $\lim_{n\to\infty} E[\rho(\Phi,\Phi_{n}^{o})] = 0$.  Next, let $\Phi_n^{*}$ be a regular conditional distribution for $\X_n$ given ${\cal F}_{m_n}$.  Then $E[\rho(\Phi,\Phi_n)] \le E[\rho(\Phi,\Phi_n^{*})]$ by (\ref{eq:rhoprp}), and  $E[\rho(\Phi,\Phi_n^{*})] = E[\rho(\Phi,\Phi_{n}^{o})]$ since the process is stationary.  So, $\lim_{n\to\infty} E[\rho(\Phi,\Phi_{n})] = 0$, as required.
	
	Thus consider the case that $m_n = 0$ and $a = 0$. From \cite{MW00} there is a martingale $M_n$ with stationary increments and a sequence $R_n$ for which $\Vert R_n\Vert/\sqrt{n} \to 0$, and $S_n = M_n+R_n\ w.p.1$ for all $n$.  Let 
$$
	\M_n(t) = {1\over\sqrt{\gamma_n}}M_{\lfloor \gamma_nt\rfloor}\quad{\rm and}\quad \R_n(t) = {1\over\sqrt{\gamma_n}}R_{\lfloor \gamma_nt\rfloor}
$$
for $0\leq t\leq b$.  Then, clearly $\X_n = \M_n + \R_n$ and $\R_n(t) \to 0$ in probability for each fixed $0 \le t \le b$.   Let $\Psi_n$ denote a regular  conditional distribution(RCD) for $\M_n$ given ${\cal F}_0$.  Then $\rho(\Phi,\Psi_n) \to 0\ w.p.1$, by the functional version of the martingale central limit theorem, applied conditionally;  see, for example, \cite{HH80}, Section 4.  From \cite{PU05} the (unconditional) distributions of $\X_n$ are tight.  So, the (unconditional) distributions of $\R_n$ are tight and, therefore, $\max_{0\le t\le b} |\R_n(t)| \to 0$ in probability.  The special case follows from the conditional version of Slutzky's theorem.

	Suppose now that $\gamma_n/m_n \to 0$ and $a < 0$.  Then, as above we may suppose $b>0$. Let $m_n^{*} = m_n + \lfloor \gamma_na\rfloor$ and let $n$ be so large that $m_n^{*} > 0$.  Then 
$$
	\X_n(t) = \X_{m_n^{*},\gamma_n}^{0,b-a}(t-a) - \X_{m_n^{*},\gamma_n}^{0,b-a}(-a) + \epsilon_n(t)
$$
for $a \le t \le b$, where 
$$
	\max_{a\le t\le b} |\epsilon_n(t)| \le 2\max_{m_n+\gamma_n a-1 \le k \le m_n+\gamma_n b+1} {|X_k|\over\sqrt{\gamma_n}}\to 0
$$ 
in probability.  So, it suffices to show that the conditional distribution of $\X_n^{*} := \X_{m_n^{*},\gamma_n}^{0,b-a}$ given ${\cal F}_0$ converges to the distribution of $\sigma\W$ in $D[0,b-a]$.  Let $\Phi^o$ be the distribution of $\sigma\W$ in $D[0,b-a]$,
$\Phi_n^{o}$ the  RCD for $\X_{0,\gamma_n}^{0,b-a}$ given ${\cal F}_0$, $\Phi_n^{*}$ a  RCD for $\X_{n}^{*}$ given ${\cal F}_0$, and $\Phi_n^{**}$ a  RCD for $\X_{n}^{*}$ given ${\cal F}_{m_n^{*}}$.  Then, as above 
$$
	E\left[\rho(\Phi^o,\Phi_n^{*})\right] \le E\left[\rho(\Phi^o,\Phi_n^{**})\right] = E\left[\rho(\Phi^o,\Phi_n^{o})\right] \to 0
$$
by (\ref{eq:rhoprp}), stationarity, and the special case.  \qed

	{\it Relation to strong mixing}.\quad  The condition (\ref{eq:MW}) may be compared with mixing conditions.  Let ${\cal G}_n = \sigma\{X_n,X_{n+1},\ldots\}$ and recall that the strong mixing coefficients are defined by
$$
	\alpha_n = \sup_{A\in{\cal F}_0,B\in{\cal G}_n} |P(A\cap B) - P(A)P(B)|.
$$
Then the condition of \cite{AH06} may be stated: for some $\epsilon > 0$,
$$
	E(X_1^{4}) < \infty\quad{\rm and}\quad \sum_{n=1}^{\infty} \alpha_n^{{1\over 2}-\epsilon} < \infty. \eqno{({\rm A9})}
$$

\begin{prop}  If (A9) holds, then (\ref{eq:MW}) holds.
\end{prop}

	{\sc Proof}.\quad First write
$$
	\Vert E(S_n|{\cal F}_0)\Vert = \sup_Y E\left[E(S_n|{\cal F}_0)Y\right] = \sup_Y E[S_nY],
$$
where the supremum is taken over all ${\cal F}_0$-measurable functions $Y$ for which $\Vert Y\Vert \le 1$.  By standard mixing inequalities, e.g., Corollary A.2 of  \cite{HH80}, Appendix III,
$$
	|E(S_nY )| \le  \sum_{k=1}^n |E(X_kY)| \le \sum_{k=1}^n 8\Vert X_k\Vert_4 \Vert Y\Vert_2 \alpha_k^{{1\over 4}} \le 8\Vert X_0\Vert_4 \sum_{k=1}^n \alpha_k^{{1\over 4}}
$$
for ${\cal F}_0$-measurable function $Y$ with $\Vert Y\Vert_2  \le 1$, where $\Vert \cdot \Vert_p$ denotes the norm in $L^p$.  So, $\Vert E(S_n|{\cal F}_0)\Vert \le 8\Vert X_0\Vert_4 \sum_{k=1}^n \alpha_k^{{1\over 4}}$ and 
$$
	\sum_{n=1}^{\infty} n^{-{3\over 2}}\Vert E(S_n|{\cal F}_0)\Vert \le 24\Vert X_0\Vert_4 \sum_{k=1}^{\infty} {1\over\sqrt{k}} \alpha_k^{{1\over 4}}.
$$
Now let $\epsilon$ be as in (A9); take $\max\{ 1,2-4\epsilon \} < q < 2$; and let $p = q/(q-1)$.  Then $p > 2$, and the right hand side of last line is at most
$$
	24\Vert X_0\Vert_4 \left[ \sum_{k=1}^{\infty} \left( {1\over {k}}\right)^{{1\over 2}p}\right]^{1\over p} \left[ \sum_{k=1}^{\infty} \alpha_k^{{1\over 4}q}\right]^{1\over q},
$$
which is finite.\qed

	Figure \ref{fig:glblac} shows the autocorrelation plot of the residual global temperature anomalies.  This is consistent with a low order autoregressive model for which (\ref{eq:MW}) is easily verified.  By way of contrast, a low order autoregressive process needs not be strongly mixing.  (The Bernoulli shift process in \cite{MW00} provides an example.) 

\begin{figure}[htbp] 
   \centering
   {\includegraphics[width=3.7in, angle=270]{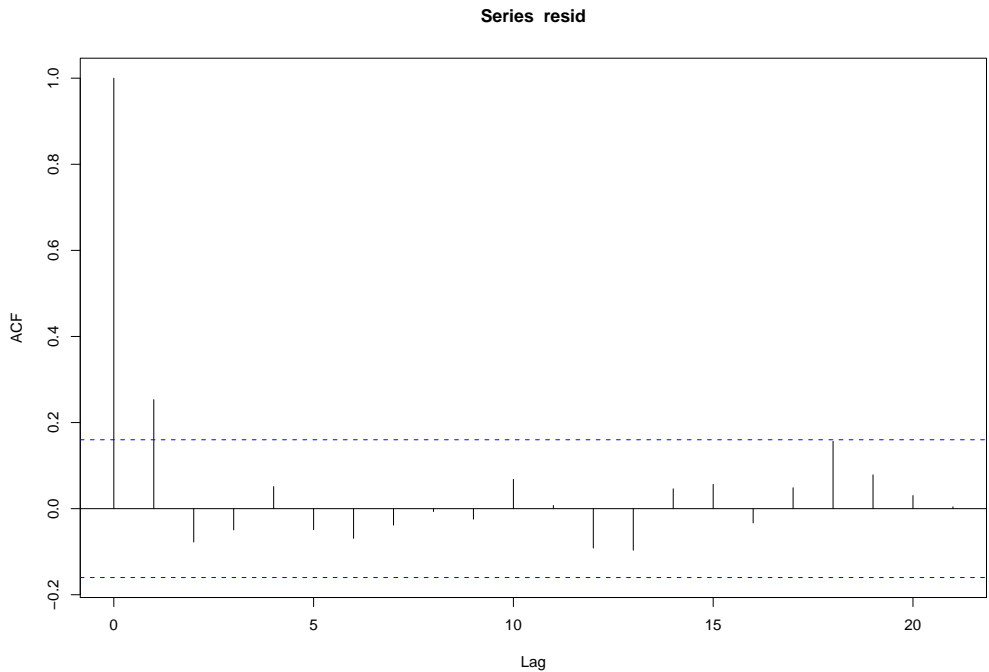}}    
   \caption{Global Temperature Anomalies}
   \label{fig:glblac}
\end{figure}

\section{Asymptotic distributions}\label{sect:ad} \hspace{1cm}

	{\it The LSE's}.\quad  Throughout this section, we suppose that the trend $\mu_k$ changes gradually in the sense that 
\be\label{eq:gc}
	\mu_k = \phi \left({k\over n}\right),
\ee
where $\phi$ is a continuous, nondecreasing function on $[0,1]$.  Thus, $\mu_k$ depends on $n$ as well as $k$, but the dependence on $n$ will be suppressed in the notation.  Let 
$$
	\Phi_n(t) = {1\over n}\sum_{j=1}^{\lfloor nt\rfloor} \phi\left({j\over n} \right),
$$
$$
	\Phi(t) = \int_{0}^t \phi(s)ds
$$
and
$$
	\tilde\phi_n(t) = \tilde{Y}_n'(t),
$$
the left hand derivative of the greatest convex minorant of $Y_n$, for $0 \le t 
\le 1$.  Then $\tilde\mu_k = \tilde\phi_n({k/n})$ and
$$
	\sup_{0\le t\le 1} |\Phi_n(t)-\Phi(t)| = O\left({1\over n} \right).
$$

A sequence $\{t_n\} \subset (0,1)$ is called {\it regular} if either $t_n \to t_0 \in (0,1)$ as $n \to \infty$ or $t_n \to 1$ and $n^{1\over 3}(1-t_n) \to \ell \in (0,\infty]$.  The first theorem obtains the asymptotic distribution of 
\be\label{eq:dlta}
	\Xi_n = n^{1\over 3}[\tilde\phi_n(t_n)-\phi(t_n)]
\ee
for regular sequences $\{ t_n \}$. Observe that if $\phi$ is continuously differentiable near $t_0$, then the 
asymptotic distribution, if any, is unchanged if $t_n$ is replaced by $[nt_n]/n$. So, for a regular sequence $\{t_n\}$, 
we implicitly assume $nt_n$ are integers with each $n=1,2,\ldots$ .

	Let $\W$ be a standard two-sided Brownian motion as in Section \ref{sect:prelim},
\be\label{eq:zee}  
		\Z(s) = \sigma\W(s) + {1\over 2}\phi'(t_0)s^2
\ee
for $s\in \R$, and 
\be\label{eq:zeen}
	\Z_n(s) = n^{2\over 3}\left[Y_n(t_n+n^{-{1\over 3}}s)-Y_n(t_n) - \phi(t_n)n^{-{1\over 3}}s\right]
\ee
for $s \in I_n := [-n^{1\over 3}t_n,n^{1\over 3}(1-t_n)]$. Then $\Xi_n = \tilde{\Z}_n'(0)$, the left hand derivative of the greatest convex minorant of $\Z_n$ at $s = 0$.

\begin{prop}\label{prop:zee}
	Suppose that (\ref{eq:MW}) and (\ref{eq:gc}) hold,  $\phi$ is continuously differentiable near $t_0 \in (0,1]$, and that $\phi'(t_0) > 0$.  Let~ $t_n \to t_0$ be regular and $0< \ell = \lim_{n\to\infty} n^{1\over 3}(1-t_n) \le \infty .$   Then  for $-\infty < a < b < \ell$, the conditional distributions of $\Z_n$ given ${\cal F}_0$ converge in $D[a,b]$ to the (unconditional) distribution of $\Z \vert [a,b]$. 
\end{prop}

	{\sc Proof}.\quad   To begin, write
\be\label{eq:zeedcmp}
	\Z_n(s) = \Psi_n(s) + \W_n(s) + R_n(s),
\ee
where
$$
	\displaylines{\Psi_n(s) = n^{2\over 3}\left[\Phi(t_n+n^{-{1\over 3}}s) - \Phi(t_n) - \phi(t_n)n^{-{1\over 3}}s\right],\cr
	\W_n(s) = n^{1\over 6}\left[\B_n(t_n+n^{-{1\over 3}}s)-\B_n(t_n)\right],\cr}
$$
and
$$
	R_n(s) = n^{2\over 3} \left[ (\Phi_n-\Phi)(t_n+n^{-{1\over 3}}s)-(\Phi_n-\Phi)(t_n) \right].
$$
It is clear that  $\sup_{s\in I_n} |R_n(s)| \le 2n^{2\over 3}\sup_{0\le t\le 1} |\Phi_n(t)-\Phi(t)| = O(n^{-{1\over 3}}) \to 0$ as $n \to \infty$ and that
$$
	\lim_{n\to\infty} \Psi_n(s) = {1\over 2}\phi'(t_0)s^2
$$
uniformly on compactas.  So, it suffices to show that the conditional distribution of $\W_n\vert [a,b]$ given ${\cal F}_0$ converges to the distribution of $\sigma \W \vert [a,b]$ in $D[a,b] $ for all compact subintervals $[a,b] \subseteq (-\infty,\ell]$; and this follows easily from Proposition \ref{prop:cc}.  To see how, let $m_n = \lfloor nt_n\rfloor$, $\gamma_n = n^{2\over 3}$, and observe that
$$
	\W_n(s) = \left[{S_{m_n+\lfloor \gamma_n s\rfloor} - S_{m_n}\over \sqrt{\gamma_n}}\right] + \epsilon_n'(s) = \X_{m_n,\gamma_n}^{a,b}(s) + \epsilon_n'(s),
$$
where $\max_{a \le s\le b} |\epsilon_n'(s)| \to 0$ in probability. \qed

	Unfortunately, $\Xi_n$ is not quite a continuous functional of $\Z_n$.  The following two lemmas are needed to obtain its limiting distribution.  The first is  simply a restatement of Lemmas 5.1 and 5.2 of \cite{WW07}.  If $f:I \to \R$ is a bounded function and $J \subseteq I$ is a subinterval, let  $G_Jf$ denote  the greatest convex minorant of $f\vert J$.

\begin{lem}\label{lem:lcl1}
	Let $f$ be a bounded piecewise continuous function on a closed interval $I$ and $[a_1,a_2] \subseteq [b_1,b_2] \subseteq I$.  If
$$
	f\left( {a_i+b_i\over 2} \right) < {G_If(a_i)+G_If(b_i)\over 2},\ i = 1,2,
$$
then $G_If = G_{[b_1,b_2]}f$ on $[a_1,a_2]$.
\end{lem}

\begin{lem}\label{lem:lcl2}
	With the notations and conditions of Proposition \ref{prop:zee},  ~ $\sup_{s\in I_n} [| \W_n(s)|  - $ $ \epsilon \min(s^2,|s|)]$ is stochastically bounded for any $\epsilon > 0$.
\end{lem}

	{\sc Proof}.\quad  Let  $I_n^{+} = [0,n^{1\over 3}(1-t_n)]$~ and~  $I_n^{-} = [-n^{1\over 3}t_n,0]$. ~ It will be shown that $\sup_{s\in I_n^{+}} |W_n(s)| - \epsilon\min(s^2,|s|) $ is stochastically bounded, the treatment of $\sup_{s\in I_n^{-}} [|W_n(s)|$ $ - \epsilon \min(s^2,|s|)]$ being similar.  Let $m_n(s) = \lfloor nt_n+n^{2\over 3}s\rfloor - \lfloor nt_n\rfloor$ for $s \in I_n^{+}$ and observe that $n^{2\over 3}s - 1 \le m_n(s) \le n^{2\over 3}s $.  Then
$$
	P\left[\sup_{s\in I_n^{+}} |\W_n(s)| - \epsilon\min(s^2,|s|) > c\right] = P\left[\sup_{s\in I_n^{+}} 
	\frac{|S_{m_n(s)}|}{n^{1\over 3}} - \epsilon\min(s^2,|s|) > c\right]
$$ for fixed $n \ge 1$ and $c > 0$.  If $c \ge 2$, the term on the right is at most 
\be\label{eq:bnd}
	 P\left[\max_{m\le n^{2\over 3}} |S_m| > n^{1\over 3}c\right] + P\left[\max_{n^{2\over 3} \le m\le n} |S_m| - n^{-{1\over 3}}\epsilon m > {1\over 2}n^{1\over 3}c\right].
\ee
Then, using the maximal inequality (\ref{eq:PU}), the right side of (\ref{eq:bnd}) is at most
\begin{align*}
	P\left[\max_{m\le n^{2\over 3}} |S_m| > {1\over 2}n^{1\over 3}c\right] + \sum_{k=1}^{\infty} &P\left[\max_{m\le 2^kn^{2\over 3}} |S_m| > {1\over 2}\left(c+\epsilon 2^k \right)n^{1\over 3}\right] \\
		&\le 24 \left[ \Vert  X_1\Vert^2+\Gamma \right] \sum_{k=0}^{\infty} {2^k\over [c+\epsilon ( 2^k-1) ]^2}, 
\end{align*}
which is independent of $n$ and approaches $0$ as $c \to \infty$. \qed

\begin{prop} \label{prop:lclz}
	If the assumptions of Proposition \ref{prop:zee} hold with $\ell=\infty$, then for any compact interval $[a_1,a_2] \subseteq \R$ and any $\epsilon > 0$, there is a compact interval $[b_1,b_2] \supseteq [a_1,a_2]$  such that
\begin{equation} \label{eq:touching}
	P\left[\tilde{\Z}_n = G_{[b_1,b_2]}\Z_n\ {\rm on}\ [a_1,a_2]\right] \ge 1-\epsilon \leq P\left[ G_{\mathbb{R}} \Z =G_{[b_1,b_2]} \Z ~~~{\rm on}~~[a_1,a_2]  \right]
\end{equation}
for all large $n$.
\end{prop}

	{\sc Proof}.\quad   Observe that $\Psi_n$ is convex in (\ref{eq:zeedcmp}) and let $\gamma = {1\over 2}\phi'(t_0)$.  Then there are $n_0 \ge 1$ and $\delta > 0$ for which  $n_0 > 1/\delta^3$ and 
$$
	{9\over 5}\gamma \le \phi'(t) \le {11\over 5}\gamma
$$
whenever $|t-t_n| \le \delta$ and $n \ge n_0$.  It then follows from a Taylor series expansion and convexity that for $n \ge n_0$
$$
	{9\over 10}\gamma s^2 \le \Psi_n(s) \le {11\over 10}\gamma s^2
$$
for $|s| \le \delta n^{1\over 3}$ and
\be\label{eq:ineq1}
	\Psi_n(s) \ge {9\over 10} \gamma \min \left( s^2,|s|\right)
\ee
for all $s \in I_n $.  Given $\epsilon$, there is a $c$ such that for all large $n$, 
\be\label{eq:ineq2}
	P\left[|\W_n(s)| + |R_n(s)|  \le {1\over 10} \gamma \min(s^2,|s|) + c\ {\rm for\ all}\ s \in I_n\right] \ge 1-\epsilon
\ee
by Lemma \ref{lem:lcl2}.  Let $B_n$ be the event defined on the left side of (\ref{eq:ineq2}).  Then $B_n$ implies $\Z_n(s) \ge 8\Psi_n(s)/9 - c$ for all $s \in I_n$ and, therefore, $\tilde{\Z}_n(s) \ge  8\Psi_n(s)/9 - c$ for all $s \in I_n$, since $\Psi_n$ is convex. Let $[a_1, a_2]$ be as in the statement of the proposition; let  $b_2 > \max \{0,a_2 \}$ be so large that 
$$
	\gamma \left[a_2^{2}+b_2^{2} - 6a_2b_2 \right] > 20c;
$$
and let $n > n_0$ be so large that $ a_2, b_2 \in I_n $ and $\max\{|a_2|,|b_2|\} \le \delta n^{1\over 3}$.   
Then $B_n$ implies
\begin{align*}
	2\Z_n\left({a_2+b_2\over 2}\right) &- \left[ \tilde{\Z}_n(a_2)+ \tilde{\Z}_n(b_2)\right] \\
			&\le  {12\over 5}\gamma\left({a_2+b_2\over 2}\right)^2 -{4\over 5} \gamma\left[a_2^{2} + b_2^{2} \right] + 4c  \\
			&= -{\gamma\over 5}\left[a_2^{2}+b_2^{2} - 6a_2b_2\right] + 4c,
\end{align*}
which is negative by the choice of $b_2$.  Similarly, for large $n$, $B_n$ implies the existence of  $b_1 < a_1$ for which $2\Z_n[{1\over 2}(a_1+b_1)] < \tilde{\Z}_n(a_1)+\tilde{Z}_n(b_1)$.  The left side of (\ref{eq:touching}) then follows from Lemma \ref{lem:lcl1}; the right hand inequality is similar, but simpler. \qed

\begin{thm}\label{thm:zeecnv}
	If the assumptions of Proposition \ref{prop:zee} hold, then the conditional distributions of $(G_{I_n}\Z_n) \vert J$ given ${\cal F}_0$ converges in probability to the distribution of $(G_{(-\infty,\ell )}\Z) \vert J $ for every compact interval $J \subseteq (-\infty,\ell )$, and the conditional distributions of  $\Xi_n$ given ${\cal F}_0$ converge in probability to the distribution of  $ [G_{(-\infty,\ell )} \Z ]'(0)$.
\end{thm}  

	{\sc Proof}.\quad  We first consider the case $\ell=\infty$.  If $J$ is any compact interval and $\epsilon > 0$, then there is a compact $K$ such that 
\be \label{eq:thmeq}
	P\left[G_{I_n}\Z_n = G_K\Z_n\ {\rm on} ~ J\right] \ge 1-\epsilon \leq P \left[ G_{\mathbb{R}}\Z=G_{K} \Z ~{\rm on}~ J \right]
\ee
for all large $n$.
Let $\Phi$ and $\Phi^{o}$ denote the distributions of $(G_{\R}\Z) \vert J$ and $(G_{K}\Z) \vert J$; and let $\Phi_n$ and $\Phi_n^{o}$ denote regular conditional distributions for $(G_{I_n}\Z_n) \vert J$ and $(G_K\Z_n) \vert J$ given ${\cal F}_0$. Recalling $\rho$ as defined in (\ref{eq:rho}) with $\mathcal{X}=D(J)$, then $E[\rho(\Phi^o,\Phi_n^{o})] \to 0$, by the continuous mapping theorem since the conditional distribution of  $\Z_{n} \vert K$ converges to the distribution of $\Z \vert K$.  It follows that
$$
	E[\rho(\Phi,\Phi_n)] \le \rho(\Phi,\Phi^o) + E[\rho(\Phi^o,\Phi_n^o)]  + E[\rho(\Phi_n^{o},\Phi_n)] \leq 3\epsilon
$$
for sufficiently large $n$, since $\rho(\Phi, \Phi^{o}) \leq \epsilon $ and  $ \rho(\Phi_n^{o}, \Phi_n)\leq \epsilon$ w.p.1 by Proposition \ref{prop:lclz}.

Now suppose $\ell<\infty$, and consider $J=[a_1,a_2]\subseteq (-\infty, \ell]$. Following the proof of Proposition \ref{prop:lclz}, there is $b_1<a_1$  for which (\ref{eq:thmeq}) holds with $K=[b_1,n^{{1\over 3}} (1-t_n)]$ for all large $n$, then the rest of the argument is similar as above. The second assertion of the theorem is an immediate consequence of the continuous mapping theorem.  \qed

\begin{cor}
Suppose (\ref{eq:MW}) and (\ref{eq:gc}) hold, then for $t\in (0,1)$,
\be\label{eq:interior}
   n^{{1\over 3}} \left( \frac{\tilde{\phi}_n(t)-\phi(t)}{\kappa} \right) \Rightarrow 2 \underset{-\infty <s<\infty}{\arg\min} \left[\W(s)+s^2 \right],
\ee
where $\kappa=[{1\over 2} \sigma^2 \phi^{\prime}(t)]^{{1\over 3}}$; moreover, for $0<\ell<\infty$,
\be\label{eq:bndry}
   n^{{1\over 3}} \left[ \tilde{\phi}_n(1-\ell n^{-{1\over 3}})-\phi(1)\right]\Rightarrow \left[ G_{(-\infty, \ell]} \left(\sigma \W(s)+{1\over 2} \phi^{\prime}(1)s^2 \right)\right]^{\prime}(0) -\ell \phi^{\prime}(1).
\ee
\end{cor}	

{\sc Proof.}\quad  The convergence follows directly from Theorem \ref{thm:zeecnv} since the left side of (\ref{eq:interior}), for example, is simply $\Xi_n/\kappa$ by taking $t_n\equiv t$.  That $ [G_{\mathbb{R}}\mathbb{Z}]^{\prime}(0) = 2\kappa  ~ {\rm argmin}_s \W(s)$ $ +s^2$ in distribution follows from rescaling properties of Brownian motion. \qed  	
\bigskip

	{\it The penalized LSE}.\quad  Now consider the penalized LSE.  Clearly,
$$
	\hat\mu_{p,n} - \mu_n = \max_{1\le k\le n} {y_{n-k+1}+\cdots+y_n -(\lambda_n+k)\mu_n\over k+\lambda_n}.
$$
The numerator here may be written  as
$$
	y_{n-k+1}+\cdots+y_n -(\lambda_n+k)\mu_n = n^{1\over 3}\left[\W_{p,n} \left(k/n^{{2\over 3}} \right) - \Delta_{n} 
\left( k/n^{{2\over 3}} \right) - n^{-{1\over 3}}\lambda_n\mu_n\right],
$$
where
$$
	\W_{p,n}(t) = n^{-{1\over 3}}\sum_{j=1}^{\lfloor n^{2\over 3}t\rfloor} X_{n-j+1}
$$
and
$$
	\Delta_n(t) = n^{-{1\over 3}}\sum_{j=1}^{\lfloor n^{2\over 3}t\rfloor} (\mu_n-\mu_{n-j+1}).
$$
It is clear that the conditional distribution of $\W_{p,n}$ converges to the distribution of $\sigma\W\ {\rm in}\ D[0,a] $ for all $0 < a < \infty$.  If (\ref{eq:gc}) holds and $\phi$ is continuously differentiable near $1$, then
$$
	\lim_{n\to\infty} \Delta_n(t) =  {1\over 2}\phi'(1)t^2
$$
uniformly on compact subintervals of $[0,\infty)$; and if $\phi'(1) > 0$, then there is an $\eta > 0$ for which $\Delta_n(t) \ge 2\eta t^2$ for all $ 6n^{-{2\over 3}} \leq t \le n^{1\over 3}$ and $ n\geq 6$.  Suppose now that $\lambda_n = \alpha n^{1\over 3}$ for some $0 < \alpha < \infty$ and let
$$
	\Z_{p,n}(t) = {\W_{p,n}(t)-\Delta_n(t) -\alpha  \phi(1) \over t+ \alpha/n^{1/3}}
$$
for $0 \le t \le n^{1\over 3}$.  Then $n^{1\over 3}(\hat\mu_{p,n}-\mu_n) = \max_{n^{-2/3} \le t\le n^{1/3}} \Z_{p,n}(t)$.

\begin{thm}
	Suppose that (\ref{eq:MW}) and (\ref{eq:gc}) hold, that $\phi$ is continuously differentiable near $1$, and that $\phi(1) \phi'(1) > 0$.  Then
$$
	n^{1\over 3}(\hat\mu_{p,n}-\mu_n) \Rightarrow \sup_{0< t<\infty} \Z_{p,\infty}(t),
$$
where
$$
	\Z_{p,\infty}(t) = {\sigma\W(t) - \alpha \phi(1) - {1\over 2}\phi'(1)t^2\over t}
$$
for $0 < t < \infty$.
\end{thm}

	{\sc Proof}.\quad  Clearly, $\Z_{p,n} \Rightarrow \Z_{p,\infty}$ in $D(K)$ for all compact subintervals $K \subseteq (0,\infty)$.  So, it suffices to show that for every $\epsilon \in (0,1)$ there is a $\delta > 0$ for which
$$
	P\left[\sup_{n^{-{2\over 3}}<t< n^{1\over 3}} \Z_{p,n}(t) = \sup_{\delta <t< \delta^{-1}} \Z_{p,n}(t)\right] \ge 1 - \epsilon
$$
for all large $n$; and for this it suffices to show that for every $\epsilon \in (0,1)$, there is a $\delta > 0$ for which
\be\label{eq:tails}
	P\left[\sup_{ n^{-2/3}<t<\delta} \Z_{p,n}(t) > -{1\over\epsilon}\right] + P\left[\sup_{\delta^{-1} \le t \le n^{1\over 3}} \Z_{p,n}(t) > -{1\over\epsilon}\right] \le \epsilon
\ee
for all large $n$.  The first term on the left side of (\ref{eq:tails}) is easy.  If $\delta < \alpha\epsilon/2$, then
$$
	P\left[\sup_{ n^{-2/3} <t<\delta} \Z_{p,n}(t) > -{1\over\epsilon}\right] \le P\left[\sup_{0 \leq t \le\delta} \W_{p,n}(t) > {1\over 2}\alpha\right],
$$
which is less than $\epsilon/2$ for all large $n$ if $\delta$ is sufficiently small, since $\W_{p,n} \Rightarrow \sigma\W$ in $D[0,1]$ and $P[\sup_{0\le t\le \delta} \W(t) > {1\over 2}\alpha] = 2[1-\Phi(\alpha/(2\sqrt{\delta})]$.  
	For the second, recall that there is an $\eta > 0$ for which  $\Delta_n(t) \ge 2\eta t^2$ for all $t \le n^{1\over 3}$ and consider $\delta < \epsilon\eta$ and $n > [\epsilon \phi(1)]^{-3}$.  Then
$$
	\eqalign{P\left[ \sup_{\delta^{-1} < t \le n^{1\over 3}} \Z_{p,n}(t) > -{1\over
\epsilon}\right] &= P\left[\W_{p,n}(t) > \Delta_n(t) - {t\over\epsilon},\ {\rm
for\ some}\ t \in [{1\over\delta},n^{1\over 3}]\right]\cr
	&\le P\left[S_n-S_{n-j} > \eta \left({j^2\over n}\right) \ {\rm
for\ some}\ j \in [{1\over\delta}n^{2\over 3},n]\right]\cr}
$$
Let $m = \lfloor \delta^{-1}n^{2\over 3}\rfloor$  Then by stationarity, the last term is at most
$$
	\sum_{k=1}^{\infty} P\left[\max_{j\le m2^k} S_j > {\eta m^2 2^{2k-2}
	\over n}\right] \le 6 \left[\Vert X_0\Vert^2 + \Gamma \right]\sum_{k=1}^{\infty} {n^2
	\over \eta^2 m^3 2^{3k-4}} \le 96 \left[ \Vert X_0\Vert^2 + \Gamma
 \right] {\delta^3\over \eta^2}
$$
for large $n$, and this may be made less than $\epsilon /2$ by taking
$\delta$ sufficiently small. \qed

\section{Simulations}\label{sect:sims}

	Simulations were conducted to assess the accuracy of the approximation implicit in (\ref{eq:interior}).  Several things affect this, including the nature of the process $\ldots X_{-1},X_0,X_1,\ldots$, the function $\phi$, the choice of $t_0$, and the sample size.   For the fluctuations, we considered an autoregressive process $X_k = \rho X_{k-1} + \epsilon_k$, where $\ldots \epsilon_{-1},\epsilon_0,\epsilon_1,\ldots$ are i.i.d. normally distributed random variables with mean $0$.  We considered two values of $\rho$, $\rho = .5$ and $.9$ representing moderate and strong dependence, three $\phi$'s, $\phi(t) = t^2,\ t$, and $\sqrt{t}$, and three values of $t_0 = 1/3,\ 1/2$, and $2/3$.  In each case the variance of $\ldots \epsilon_{-1},\epsilon_0,\epsilon_1,\ldots$ was chosen to make $\sqrt{E[X_k^2]}=.25$.  The sample size was $n = 150$.  In many ways, these choices are consistent with the global warming example.  For each choice of these values $10,000$ time series were generated and the empirical distribution function of $\Xi_n/\kappa$ was computed at selected percentiles of Chernoff's distribution \cite{GW01}.  The results are presented in the Tables 1, 2, and 3.

\vskip.5in

{\def\baselinestretch{1.00}
\begin{center}
\begin{table}[h]
\label{tbl:square}
\caption{$\phi(t) = t^2$}
\begin{tabular}{|cc||ccccc||ccccc|}

\hline
&&&&&&&&&&& \\
&&&&$\rho=.5$&&&&&$\rho=.9$&& \\
&&&$t_0$&&&&&$t_0$&&&\\

$p$ & $\pm$ & $1/3$ & $1/2$ & $2/3$ & min  & max & $1/3$& $1/2$ & $2/3$ & min  & max \\
\hline

  .025 & .0016 & .0142 & .0181 & .0225 & .0142 & .0225 & .0219 & .0104 & .0119 & .0102 & .0219 \\
  .050 & .0022 & .0341 & .0406 & .0448 & .0323 & .0448 & .0443 & .0275 & .0303 & .0262 & .0443 \\
  .100 & .0030 & .0773 & .0898 & .0939 & .0762 & .0950 & .0878 & .0676 & .0734 & .0673 & .0878 \\
  .200 & .0040 & .1725 & .1907 & .1892 & .1717 & .1920 & .1808 & .1631 & .1660 & .1584 & .1808 \\
  .250 & .0043 & .2231 & .2373 & .2398 & .2219 & .2440 & .2271 & .2140 & .2193 & .2067 & .2271 \\
  .300 & .0046 & .2745 & .2906 & .2892 & .2706 & .2951 & .2761 & .2643 & .2723 & .2563 & .2761 \\
  .400 & .0049 & .3781 & .3864 & .3902 & .3762 & .3962 & .3719 & .3683 & .3827 & .3581 & .3827 \\
  .500 & .0050 & .4834 & .4915 & .4930 & .4786 & .4963 & .4738 & .4718 & .4925 & .4673 & .4929 \\
  .600 & .0049 & .5908 & .5921 & .5962 & .5826 & .5990 & .5769 & .5805 & .5955 & .5729 & .6002 \\
  .700 & .0046 & .6957 & .6893 & .6986 & .6852 & .6997 & .6765 & .6847 & .7068 & .6765 & .7068 \\
  .750 & .0043 & .7448 & .7381 & .7470 & .7381 & .7521 & .7271 & .7400 & .7579 & .7271 & .7591 \\
  .800 & .0040 & .7943 & .7942 & .7977 & .7889 & .8021 & .7818 & .7924 & .8087 & .7811 & .8121 \\
  .900 & .0030 & .8947 & .8985 & .9014 & .8888 & .9021 & .8863 & .8980 & .9050 & .8863 & .9093 \\
  .950 & .0022 & .9455 & .9491 & .9508 & .9442 & .9532 & .9383 & .9502 & .9556 & .9382 & .9556 \\
  .975 & .0016 & .9709 & .9746 & .9763 & .9709 & .9781 & .9663 & .9759 & .9801 & .9660 & .9801 \\
  \hline
  \end{tabular}
  
 \vskip.1in
{\bf Note}: Columns three, four, and five show the empirical distribution function of scaled $\Xi_n$ at the $p^{\rm th}$ percentile of Chernoff's distribution for $t_0 = 1/3,\ 1/2$, and $2/3$.  The value of $p$ is in column one, and column two lists the standard errors of the simulations.  Columns six and seven list the minimum and maximum of the empirical distribution function over $1/3 \le t_0 \le 2/3$.  Columns eight through twelve provide the same information for $\rho = .9$

\end{table}

\vskip.25in

 \begin{table}
 \label{tbl:lin}
 \caption{$\phi(t) = t$}
 \begin{tabular}{|cc||ccccc||ccccc|}
\hline
&&&&&&&&&&& \\
&&&&$\rho=.5$&&&&&$\rho=.9$&& \\
&&&$t_0$ &&&&&$t_0$&&& \\
$p$ & $\pm$ & $1/3$ & $1/2$ & $2/3$ & min  & max & $1/3$& $1/2$ & $2/3$ & min  & max \\
\hline
   .025 & .0016 & .0417 & .0222 & .0135 & .0129 & .0417 & .0394 & .0202 & .0107 & .0104 & .0394 \\
  .050 & .0022 & .0709 & .0468 & .0324 & .0320 & .0709 & .0674 & .0417 & .0266 & .0266 & .0674 \\
 .100 & .0030 & .1246 & .0946 & .0731 & .0731 & .1246 & .1223 & .0880 & .0652 & .0652 & .1223 \\
  .200 & .0040 & .2226 & .1970 & .1713 & .1713 & .2227 & .2185 & .1834 & .1601 & .1601 & .2185 \\
 .250 & .0043 & .2729 & .2472 & .2214 & .2214 & .2729 & .2642 & .2328 & .2105 & .2105 & .2651 \\
 .300 & .0046 & .3169 & .2971 & .2723 & .2723 & .3169 & .3099 & .2848 & .2618 & .2618 & .3105 \\
 .400 & .0049 & .4048 & .3987 & .3808 & .3808 & .4051 & .4044 & .3923 & .3781 & .3736 & .4044 \\
 .500 & .0050 & .4923 & .5023 & .4899 & .4878 & .5058 & .4925 & .4983 & .4955 & .4877 & .5014 \\
 .600 & .0049 & .5852 & .6050 & .6036 & .5820 & .6108 & .5851 & .6010 & .6111 & .5851 & .6154 \\
 .700 & .0046 & .6785 & .7083 & .7150 & .6726 & .7194 & .6786 & .7052 & .7242 & .6786 & .7253 \\
 .750 & .0043 & .7233 & .7567 & .7694 & .7194 & .7733 & .7254 & .7612 & .7826 & .7254 & .7826 \\
 .800 & .0040 & .7718 & .8078 & .8215 & .7644 & .8249 & .7735 & .8107 & .8352 & .7735 & .8367 \\
 .900 & .0030 & .8666 & .9031 & .9222 & .8666 & .9225 & .8737 & .9087 & .9295 & .8737 & .9304 \\
 .950 & .0022 & .9245 & .9523 & .9667 & .9245 & .9689 & .9302 & .9564 & .9716 & .9302 & .9716 \\
 .975 & .0016 & .9562 & .9755 & .9877 & .9562 & .9877 & .9615 & .9800 & .9887 & .9615 & .9887  \\
 \hline
 \end{tabular}
 \vskip.1in
 {\bf Note}: See note to Table 1.
\end{table}

\vskip.25in

\begin{table}
\label{tbl:sqrt}
\caption{$\phi(t) = \sqrt{t}$}
\begin{tabular}{|cc||ccccc||ccccc|}
\hline
&&&&&&&&&&& \\
&&&&$\rho=.5$&&&&&$\rho=.9$&& \\
&&&$t_0$&&&&&$t_0$&&& \\
$p$ & $\pm$ & $1/3$ & $1/2$ & $2/3$ & min  & max & $1/3$& $1/2$ & $2/3$ & min  & max \\
\hline

 .025 & .0016 & .0392 & .0130 & .0046 & .0046 & .0392 & .0444 & .0155 & .0040 & .0040 & .0444 \\
 .050 & .0022 & .0695 & .0309 & .0151 & .0146 & .0695 & .0758 & .0328 & .0138 & .0138 & .0758 \\
 .100 & .0030 & .1259 & .0770 & .0456 & .0456 & .1259 & .1292 & .0731 & .0422 & .0422 & .1292 \\
 .200 & .0040 & .2270 & .1759 & .1354 & .1346 & .2270 & .2291 & .1742 & .1371 & .1371 & .2291 \\
 .250 & .0043 & .2709 & .2284 & .1909 & .1892 & .2716 & .2771 & .2209 & .1907 & .1907 & .2780 \\
 .300 & .0046 & .3193 & .2840 & .2496 & .2478 & .3219 & .3288 & .2819 & .2484 & .2484 & .3288 \\
 .400 & .0049 & .4176 & .3925 & .3772 & .3759 & .4211 & .4263 & .3940 & .3723 & .3723 & .4271 \\
 .500 & .0050 & .5088 & .5037 & .5060 & .5021 & .5194 & .5216 & .5146 & .5039 & .5028 & .5217 \\
 .600 & .0049 & .6071 & .6154 & .6375 & .6050 & .6375 & .6111 & .6261 & .6296 & .6111 & .6343 \\
 .700 & .0046 & .6979 & .7291 & .7581 & .6979 & .7620 & .7048 & .7386 & .7633 & .7048 & .7633 \\
.750 & .0043 & .7478 & .7831 & .8180 & .7478 & .8180 & .7532 & .7890 & .8198 & .7532 & .8198 \\
.800 & .0040 & .7966 & .8334 & .8730 & .7966 & .8734 & .8052 & .8421 & .8717 & .8052 & .8717 \\
.900 & .0030 & .8881 & .9312 & .9595 & .8881 & .9599 & .9029 & .9360 & .9584 & .9029 & .9584 \\
.950 & .0022 & .9424 & .9741 & .9869 & .9424 & .9875 & .9480 & .9759 & .9839 & .9480 & .9850 \\
.975 & .0016 & .9686 & .9897 & .9966 & .9686 & .9970 & .9728 & .9909 & .9948 & .9728 & .9948  \\
\hline
 
\end{tabular}
 \vskip.1in
 {\bf Note}: See note to Table 1.
\end{table}

\vskip.25in

\begin{table}
\label{tbl:bndry}
\caption{The Normalized Boundary Corrected Estimator}
\begin{tabular}{|cc||ccc|ccc|ccc|}

\hline
&&&&&&&&&&\\
&&&$\phi=\sqrt{t}$&&&$\phi=t$&&&$\phi=t^2$ & \\
\hline
$z$ & $n=\infty$ & $\rho = 0$ &  $\rho = .5$ & $\rho = .9$ & $\rho = 0$ &  $\rho = .5$ & $\rho = .9$ & $\rho = 0$ &  $\rho = .5$ & $\rho = .9$ \\
\hline
-2.5 & .0344 & .0426 & .0447 & .0503 & .0480 & .0370 & .0299 & .0399 & .0363 & .0216 \\

-2.0 & .0936 & .1054 & .1079 & .1144 & .1144 & .0981 & .0923 & .1033 & .1025 & .0786 \\ 

-1.5 & .1913 & .2129 & .2134 & .2212 & .2358 & .2076 & .2034 & .2135 & .2087 & .1905 \\ 

-1.0 & .3242 & .3603 & .3542 & .3697 & .3846 & .3508 & .3643 & .3588 & .3625 & .3609 \\  

-0.5 & .4786 & .5137 & .5091 & .5237 & .5495 & .5135 & .5485 & .5221 & .5374 & .5551 \\ 

\phantom{-}0.0 & .6293 & .6609 & .6608 & .6878  & .6989 & .6744 & .7103 & .6769 & .6926 & .7383 \\ 

\phantom{-}0.5 & .7549 & .7784 & .7829 & .8109 & .8177 & .7972 & .8375 & .8019 & .8238 & .8724 \\  

\phantom{-}1.0 & .8483 & .8687 & .8734 & .9048 & .8988 & .8840 & .9170 & .8861 & .9047 & .9435 \\  

\phantom{-}1.5 & .9137 & .9303 & .9310 & .9570 & .9463 & .9391 & .9662 & .9390 & .9536 & .9828 \\ 

\phantom{-}2.0 & .9568 & .9647 & .9660 & .9816 & .9746 & .9692 & .9877 & .9703 & .9811 & .9944 \\ 

\phantom{-}2.5 & .9787 & .9832 & .9860 & .9935 & .9898 & .9872 & .9970 & .9874 & .9925 & .9988 \\ 

\hline

\end{tabular}
\vskip.1in
{\bf Note}: Column 2 lists a Monte Carlo estimate of the asymptotic distribution function; columns 3, 4, and 5, list estimates of the actual distribution function for $n = 150$, $\rho = 0, .5$, and $.9$, and $\phi(t) = \sqrt{t}$; columns 6, 7, and 8 provide the same information for $\phi(t)= t$, columns 9, 10, and 11 for $\phi(t) = t^2$. 
\end{table}

\vskip.25in
\begin{table}
\label{tbl:pen}
\caption{The Normalized Penalized Estimator}
\begin{tabular}{|cc||ccc|ccc|ccc|}

\hline
&&&&&&&&&&\\
&&&$\phi=\sqrt{t}$&&&$\phi=t$&&&$\phi=t^2$ & \\
\hline
$z$ & $n=\infty$ & $\rho = 0$ &  $\rho = .5$ & $\rho = .9$ & $\rho = 0$ &  $\rho = .5$ & $\rho = .9$ & $\rho = 0$ &  $\rho = .5$ & $\rho = .9$ \\
\hline
-2.5 & .0393 & .0412 & .0470 & .0577    & .0379 & .0389 & .0363      & .0340 & .0314 & .0241 \\

-2.0 & .1122 & .1153 & .1331 & .1514 & .1065 & .1176 & .1283          & .1019 & .1067 & .1037 \\ 

-1.5 & .2380 & .2424 & .2779 & .3080 &   .2354 & .2627 & .3045        & .2319 & .2596 & .2965 \\ 

-1.0 & .3952 & .4082 & .4552 & .5117 &   .4077 & .4574 & .5511 &      .4102 & .4646 & .4779 \\  

-0.5 & .5582 & .5721 & .6416 & .6637 &     .5934 & .6527 & .7592 &    .6029 & .6717 & .7085 \\ 

\phantom{-}0.0 & .7055 & .7265 & .7679 & .7985     & .7470 & .8038 & .8977 &   .7611 & .8317 & .9304 \\ 

\phantom{-}0.5 & .8169 & .8426 & .8903 & .9354     & .8553 & .9014 & .9674 &   .8674 & .9277 & .9836 \\  

\phantom{-}1.0 & .8948 & .9186 & .9509 & .9861     & .9234 & .9581 & .9916 &   .9366 & .9720 & .9975 \\  

\phantom{-}1.5 & .9401 & .9602 & .9810 & .9967     & .9650 & .9831 & .9995 &   .9718 & .9908 & .9997 \\ 

\phantom{-}2.0 & .9677 & .9821 & .9922 & .9993     & .9842 & .9932 & 1.000 &   .9884 & .9978 & 1.000 \\ 

\phantom{-}2.5 & .9837 & .9916 & .9974 & .9999     & .9923 & .9983 & 1.000 &   .9960 & .9996 & 1.000 \\ 

\hline

\end{tabular}
\vskip.1in
{\bf Note}: Column 2 lists a Monte Carlo estimate of the asymptotic distribution function; columns 3, 4, and 5, list estimates of the actual distribution function for $n = 150$, $\rho = 0, .5$, and $.9$, and $\phi(t) = \sqrt{t}$; columns 6, 7, and 8 provide the same information for $\phi(t)=  t$, columns 9, 10, and 11 for $\phi(t) = t^2$.

\end{table}

\end{center}

	In Table 1, the agreement between the empirical distribution function and the limiting distribution seems generally better in the right tail than the left where  the empirical is consistently less than the limiting distribution.  In Table 2, the agreement is excellent at $t_0 = 1/2$ but deteriorates markedly for $t_0 = 1/3$ or $2/3$.  In Table 3, the empirical distribution of the absolute value appears to be stochastically smaller 
than the corresponding limit in all but two columns ($t_0 = 1/3$).  In all three tables the empirical distribution function is generally decreasing in $t_0$.  This is easily explained by the numbers of maxima and minima in the Max-Min formula, (\ref{eq:minmax}).  Also, in all tables the difference between moderate and strong dependence is modest, suggesting that the effect of dependence is adequately captured in the calculation of $\sigma$.
	
	Similar simulations showed that the approximations implicit in Theorem 2 and (\ref{eq:bndry}) were not so good (depending on $\ell$) in the case of (\ref{eq:bndry}).  Monte Carlo estimates of the distribution function of (\ref{eq:bndry}) are listed in Table 4 for $n = 150$, $\rho = 0, .5$, and $.9$, and the same three functions $\phi$ along with the asymptotic distribution function.   Similar results were obtained for the normalized penalized estimator and are presented in Table 5.   While the agreement leaves much to be desired, the results are not without practical implications:  at the very least they suggest that the limiting distributions are not highly sensitive to the distribution of the fluctuations, within broad limits; and this suggestion is confirmed in Tables 4 and 5 which show good agreement for the three values of $\rho$.

\end{document}